\documentclass[11pt]{article}
\usepackage{amssymb,amsmath,amsthm,fullpage}

 \newtheorem{thm}{Theorem}
 \newtheorem{cor}[thm]{Corollary}
 \newtheorem{lem}[thm]{Lemma}
 \newtheorem{prop}[thm]{Proposition}
 \newtheorem{conj}[thm]{Conjecture}
 \theoremstyle{definition}
 \newtheorem{defn}[thm]{Definition}
 \theoremstyle{remark}

\begin{document}

\title{New upper bounds on the chromatic number of a graph}

\author{ landon rabern \\
         rabernsorkin@cox.net }

\date{\today}

\maketitle

\begin{abstract}
We outline some ongoing work related to a conjecture of Reed \cite{reed97} on $\omega$, $\Delta$, and $\chi$.  We conjecture
that the complement of a counterexample $G$ to Reed's conjecture has connectivity on the order of $\log(|G|)$.  We prove that this holds
for a family (parameterized by $\epsilon > 0$) of relaxed bounds; the $\epsilon = 0$ limit of which is Reed's upper bound.
\end{abstract}

\section{Some Bounds}

In all that follows, \emph{graph} will mean finite simple graph with non-empty vertex set.

We will need the following result from \cite{rabern}.

\begin{prop}
Let $I_1, \ldots, I_m$ be disjoint independent sets in a graph G. Then
\begin{equation}\label{indepdendentsetinequality} \chi (G) \leq \frac{1}{2} \left (\omega (G) + |G| - \displaystyle \sum_{j=1}^m |I_j| + 2m -1 \right). \end{equation}
\end{prop}

In \cite{rabernReed}, an upper bound on the chromatic number in terms of the clique number, maximal degree and order was given.
\begin{prop}
Let $G$ be a graph.  Then,
\[\chi(G) \leq \frac{1}{2} \left (\omega(G) + \frac{|G| + \Delta(G) + 1}{2} \right).\]
\end{prop}

We can do better than this when $G$ has an induced subgraph with order much bigger than chromatic number.

\begin{prop}
Let $G$ be a graph. Then, for any induced subgraph $H$ of $G$,
\[\chi(G) \leq \frac{1}{2} \left (\omega(G) + \frac{\Delta(G) + 1 + |G|}{2} \right ) + \frac{3\chi(H) - |H|}{4}.\]
\end{prop}

\begin{prop}
Let $G$ be a graph and $K$ a cut-set in $\overline{G}$. Then, for any induced subgraph $H$ of $G$,
\[\chi(G) \leq \frac{1}{2} \left (\omega(G) + \Delta(G) + 1 \right ) + \frac{4\chi(G[K]) + 3\chi(H \smallsetminus K) - |H \smallsetminus K|}{4}.\]
\end{prop}

\begin{cor}
Let $G$ be a graph. Then, for any induced subgraph $H$ of $G$,
\[\chi(G) \leq \frac{1}{2} \left (\omega(G) + \Delta(G) + 1 \right ) + \frac{5\kappa(\overline{G}) + 3\chi(H) - |H|}{4}.\]
\end{cor}

\begin{cor}
Let $G$ be a graph and $K$ a cut-set in $\overline{G}$. Then
\[\chi(G) \leq \frac{1}{2} \left (\omega(G) + \Delta(G) + 1 \right ) + \frac{4\chi(G[K]) + \alpha(G[K]) + 3 - \alpha(G)}{4}.\]
\end{cor}

\begin{cor}
Let $G$ be a graph. Then
\[\chi(G) \leq \frac{1}{2} \left (\omega(G) + \Delta(G) + 1 \right ) + \kappa(\overline{G}) + 1 - \frac{\alpha(G)}{4}.\]
\end{cor}

\section{Chromatic Excess}
\begin{defn}
Let $G$ be a graph.  The \emph{chromatic excess} of $G$ is defined to be
\[\eta(G) = \max_{H \leq G} |H| - 3\chi(H).\]
\end{defn}

The extremal cases of $H$ being a maximal independent set and $H = G$ give the following constraints.
\begin{lem}
$\alpha - 3 \leq \eta \leq \frac{\alpha - 3}{\alpha}n$.
\end{lem}

\begin{lem}
$\eta \geq n - 3\chi$.
\end{lem}

We may rewrite Proposition 3 and Corollary 5 in terms of the chromatic excess.

\begin{cor}
\[\chi \leq \frac{1}{2} \left (\omega + \frac{\Delta + 1 + n}{2} \right ) - \frac{\eta}{4}.\]
\end{cor}

Or equivalently.

\begin{cor}
\[\chi \leq \frac{1}{2} (\omega + \Delta + 1) + \frac{\overline{\delta} - \eta}{4}.\]
\end{cor}

\begin{cor}
\[\chi \leq \frac{1}{2} \left (\omega + \Delta + 1 \right ) + \frac{5\overline{\kappa} - \eta}{4}.\]
\end{cor}

\section{A Counterexample To Reed's Conjecture Has Highly Connected Complement?}

\begin{prop}
For every $\epsilon > 0$ there exists a constant $C(\epsilon) > 0$ such that for any graph $G$, at least one of the following holds
\begin{itemize}
\item $\chi(G) \leq (\frac{1}{2} + \epsilon)\omega(G) + \frac{\Delta(G) + 2}{2}$,
\item $\kappa(\overline{G}) \geq C(\epsilon)\log(|G|)$.
\end{itemize}
\end{prop}
\begin{proof}
Assume the former does not hold.  Apply Corollary 7 to get upper bounds on $\alpha$ and $\omega$ in terms of $\overline{\kappa}$.  Now use Ramsey Theory.
\end{proof}

\begin{conj}
There exists a constant $C > 0$ such that 
\[\chi > \left \lceil \textstyle \frac{1}{2}(\omega + \Delta + 1) \right \rceil \Rightarrow \overline{\kappa} \geq C\log(n).\]
\end{conj}

\end{document}